\documentclass[letterpaper, 10 pt, conference]{ieeeconf}  

\IEEEoverridecommandlockouts
\usepackage{graphics} 
\usepackage{epsfig} 
\usepackage{times} 
\usepackage{bbm}
\usepackage{xcolor}
\usepackage{soul}
\usepackage{float}
\usepackage{amsmath,amssymb,amsfonts}

\newtheorem{definition}{Definition}
\newtheorem{theorem}{Theorem}
\newtheorem{lemma}{Lemma}
\newtheorem{assumption}{Assumption}
\newtheorem{corollary}{Corollary}
\newtheorem{proposition}{Proposition}

\title{\LARGE \bf
Sample Complexity of Chance Constrained Optimization in Dynamic Environment\\
}

\author{Apurv Shukla, Qian Zhang and Le Xie
\thanks{This material is based upon work supported by the U.S. Department of
Energy's Office of Energy Efficiency and Renewable Energy (EERE) under the Solar Energy
Technologies Office Award Number DE-EE0009031.}
\thanks{AS, QZ, LX are with the ECE department at Texas A\&M University, College Station, TX e-mail\{apurv.shukla,zhangqianleo,le.xie\}@tamu.edu}%
}

\begin{document}

\maketitle
\thispagestyle{empty}
\pagestyle{empty}

\begin{abstract}

We study the scenario approach for solving chance-constrained optimization in time-coupled dynamic environments. Scenario generation methods approximate the true feasible region from scenarios generated independently and identically from the actual distribution. In this paper, we consider this problem in a dynamic environment, where the scenarios are assumed to be drawn in a sequential fashion from an unknown and time-varying distribution. Such dynamic environments are driven by changing environmental conditions that could be found in many real-world applications such as energy systems. We couple the time-varying distributions using the Wasserstein metric between the sequence of scenario-generating distributions and the actual chance-constrained distribution. Our main results are bounds on the number of samples essential for ensuring the ex-post risk in chance-constrained optimization problems when the underlying feasible set is convex or non-convex. Finally, our results are illustrated on multiple numerical experiments for both types of feasible sets.

\end{abstract}

\section{Introduction}
\label{sec:introduction}
Chance-constrained optimization (CCO) is a general framework for decision-making in uncertain environments when the goal is to explicitly model the risk of constraint satisfaction. Initial applications focused on production scheduling while making distributional assumptions on the family of distributions, such as multivariate Gaussian distribution \cite{charnes1958cost,geng2019data}. \\ 
\indent Over the past decade, due to the rapid deployment of sensing, communication, and computational capabilities in extremely reliable and cost-effective manners, many novel data-driven methods were proposed to solve CCO, which can be briefly classified into three categories in \cite{geng2021computing}: (1) safe approximation \cite{nemirovski2012safe,bertsimas2021probabilistic}; (2) sample average approximation (SAA) \cite{luedtke2008sample,luedtke2010integer}; and (3) the scenario approach \cite{calafiore2006scenario,campi2008exact,campi2018general}. The key difference between the safe approximation method and robust optimization is that the solution by safe approximation is \textit{guaranteed} to meet the risk level \cite{ben2009robust}. In SAA, the true distribution of the uncertain variables is approximated by the empirical distribution from collected samples. The sample average approximation allows the violation of constraints under some sampled scenarios. This violation program can be formulated as a mixed-integer program~\cite{luedtke2010integer}. An alternative popular approach for CCO is the scenario approach, which doesn't make assumptions on the underlying distribution and seeks the optimal solutions that are feasible for \textit{all} the extracted scenarios.\\
\indent The above three methods assume that the scenarios are sampled \textit{independently from identical} distributions. In practical settings, this assumption does not necessarily hold true. Over a large time horizon, the scenario-generating distributions change. For example, due to the high penetration of renewable energy, there is an increase in uncertainty during economic dispatch decisions in a power system~\cite{gu2016stochastic}. Improved forecast techniques reduce this uncertainty~\cite{xie2013short}, but the forecast uncertainty rises significantly as the forecast horizon increases, implying a changing distribution of uncertain generation. Similar examples can be found in many other fields, such as adversarial machine learning \cite{ortiz2021optimism,bienstock2022robust}.\\
\indent Recent work has focused on a robust SAA scheme for solving CCO in dynamic environments \cite{yan2022data}. This scheme inherits disadvantages from traditional SAA, including non-convex feasible regions which may be difficult to solve \cite{luedtke2008sample}. This shortage arises because sampled constraints are to be violated, and need to add more assumptions to gain the exact feasible region. Different from SAA, the scenario approach only proves the property of the optimal solution, not a region, which can easily be solved in convex problems. After first proposed in 2006 \cite{calafiore2006scenario}, fruitful breakthroughs were witnessed in scenario approach these years, including exact feasibility \cite{campi2008exact}, scenario theory for nonconvex optimization \cite{campi2018general}, and computation algorithm for essential sets \cite{geng2021computing}, but the problem with time-varying distributions is still unexplored as pointed out in~\cite{campi2021scenario}. 
\subsection{Summary of Contributions}
In this paper, we analyze the generation of scenarios from time-varying distributions. Our suggested contribution is establishing bounds on the violation probability of an optimal solution of the robust scenario problem for guaranteeing prescribed risk levels in chance-constrained optimization when the scenarios are generated with underlying time-varying distributions for both convex and non-convex feasible regions. The results for convex and non-convex optimization are of a fundamentally different nature compared to the convex case. Our guarantees are tight in the sense that we are able to recover classical bounds for these problems when scenarios are sampled from a fixed distribution as in~\cite{campi2018general} and~\cite{yan2022data}.
In this paper, we consider a model for generating scenarios drawn from distributions centered around a central distribution and coupled across time using a suitable Wasserstein ball with different radii. We build on their work to further consider the case when the data-generating distributions changes with time. Results are further elaborated on multiple numerical experiments in convex and non-convex problems. These experiments demonstrate the tightness of our results by perfectly interpolating between static~\cite{campi2018general} and dynamic environments. We show that our theoretical guarantees match the best-known results for a static environment and scale gracefully with parameters measuring the non-stationarity of the scenario-generating distribution.

The remainder of this article is organized as follows. In Section \ref{sec:scenario-generation}, we formalize the CCO in a dynamic environment. Section~\ref{sec:nonstationary-scenario} establishes guarantees for CCO in dynamic environments and numerical results are presented in Section~\ref{sec:case-study}. The notations in this article are standard. Given a real number $x \in \textbf{R}$, we denote $(x)_+:=\text{max}\{0,x\}$. The cardinality of a set $\mathcal{N}$ is $|\mathcal{N}|$. Removal of subset $I$ from set $\mathcal{I}$ is represented by $\{\mathcal{I}\setminus I\}$. The essential supremum is ess sup.

\section{Chance Constrained Optimization with Observation Robustness}
\label{sec:scenario-generation}
\subsection{Chance-Constrained Optimization (CCO)}
A typical formulation of Chance-Constrained Optimization (CCO) is~\eqref{CCO}
\begin{subequations}
\label{CCO}
\begin{align}
\min _x &\ c^T x \nonumber \\
\text { s.t. } & \mathbb{P}_{\xi}(f(x, \xi) \leq 0) \geq 1-\epsilon \nonumber \\
& g(x) \leq 0 \nonumber
\end{align}
\end{subequations}
where $\xi \in \Delta$ is a random vector defined in uncertainty set $\Delta$. We could write~\eqref{CCO} more compactly by defining $\mathcal{X}_{\xi}:=\left\{x \in \mathbf{R}^n: f(x, \xi) \leq 0\right\}$ and $\chi:=\left\{x \in \mathbf{R}^n: g(x) \leq 0\right\}$ as:
\begin{equation}
\label{CCO2}
\begin{aligned}
\min _{x \in \chi} &\ c^{\top} x \\
\text { s.t. } & \mathbb{P}_{\xi}\left(x \in \mathcal{X}_{\xi}\right) \geq 1-\epsilon
\end{aligned}
\end{equation}

\indent The scenario approach randomly extracts \textit{N independent and identically distributed (i.i.d.)} scenarios $\mathcal{S}:=\left\{\xi_1, \xi_2, \cdots, \xi_N\right\}$ to approximate the chance-constrained program (\ref{CCO}) with the following \textit{scenario problem}:
\begin{equation}
\label{SP}
\begin{aligned}
\mathrm{SP}(\mathcal{S}): & \min _x c^{\top} x\\
\text { s.t. } & f\left(x, \xi_i\right) \leq 0 &\; i=1,2,3,...,N\\
& g(x) \leq 0
\end{aligned}
\end{equation}
The scenario problem ${SP}(\mathcal{S})$ can be written by defining $\mathcal{X}_{\xi_i}:=\left\{x \in \mathbf{R}^n: f(x, {\xi}_i) \leq 0\right\}$.
\begin{equation}
\begin{aligned}
\mathrm{SP}(\mathcal{S}): &\min _{x \in \chi} c^{\top} x\\
\text { s.t. } &x \in \cap_{i=1}^N \mathcal{X}_{\xi_i}
\end{aligned}
\end{equation}

When scenarios are generated \textit{i.i.d.}, the measurement errors are affected by the inherent uncertainty due to noise by using techniques such as state estimation~\cite{xie2012fully}, Kalman Filtering~\cite{welch1995introduction}, etc.
and it is important for the scenario solution to be robust to the underlying observation error or equivalently to the observed scenarios. Formally,

we can obtain a robust scenario solution as the solution to the following min-max robust optimization problem
called \textit{robust scenario problem}, $\mathrm{RSP}(\mathcal{S})$.
\begin{equation}
\label{RSP2}
\begin{aligned}
\mathrm{RSP}(\mathcal{S}):\min _x c^{\top} x\\
\text{ s.t. } \max _{\xi_i} f\left(x, \xi_i\right) \leq 0 &&\\
\Vert {\xi}_i-{\eta}_i\Vert &\leq& r_i,\ i=1,2,3,...,N\\
g(x) \leq 0 &&
\end{aligned}
\end{equation}
Optimization problem~(\ref{RSP2}) can be equivalently reformulated by defining a compact set: 
$\mathcal{X}_{\eta_i}^{r_i}:=\left\{x \in \mathbf{R}^n: \max _{\xi_i} f\left(x, \xi_i\right) \leq 0, \; \|{\xi}_i-{\eta}_i\|\leq r_i\right\}$, as follows:
\begin{equation}
\label{eqn:robust-scenario-compact}
\begin{aligned}
\mathrm{RSP}(\mathcal{S}): &\min _{x \in \chi} c^{\top} x\\
\text { s.t. } &x \in \cap_{i=1}^N \mathcal{X}_{\eta_i}^{r_i}
\end{aligned}
\end{equation}
\indent The \textit{robust scenario approach} can be regarded as a two-stage optimization problem. The first stage consists of $N$ independent robust optimization problems with the objective function $f\left(x, \xi_i\right)$ and decision variable $\xi_i$ respectively. This stage can be efficiently solved by parallelization. The optimal solution from the first stage programs will construct the constraints of the second stage optimization, i.e., the traditional scenario approach scheme (\ref{SP}). Since $\mathcal{X}_{\xi_i} \subseteq \mathcal{X}_{\eta_i}^{r_i}$, we have that:
\begin{equation}
\mathbb{P}^N\left(\mathbb{V}\left(x_{\mathcal{RS}}^*\right)>\epsilon\right) \leq \mathbb{P}^N\left(\mathbb{V}\left(x_{\mathcal{N}}^*\right)>\epsilon\right)
\end{equation}
In $\mathrm{RSP}$, every scenario may have different measure $\mathbb{Q}_i$. Closely related Ambiguous CCO problem assumes that the scenarios are sampled from a fixed but unknown measure $\mathbb{Q}$ \cite{erdougan2006ambiguous,tseng2016random}.

\indent
We recollect some crucial definitions based on~\cite{campi2008exact,calafiore2010random} that will be of importance in the results to follow. Let $x_\mathcal{N}^*$ and $x_\mathcal{RS}^*$ stand for the optimal solution to scenario problems $\mathrm{SP}(\mathcal{S})$ and robust scenario problems $\mathrm{RSP}(\mathcal{S})$. 
\begin{definition}[Violation Probability] 
\label{def:violation-prob}
The \textit{violation probability} of a candidate solution $x^*$ is defined as the probability that $x^*$ is infeasible, i.e., $\mathbb{V}(x^*):=\mathbb{P}_\xi(x^*\notin \mathcal{X}_\xi)$ 
\end{definition}

\begin{definition}[Support Scenario] 
Scenario $\xi_i$ is a \textit{support scenario} for the scenario problem $\mathrm{SP}(\mathcal{S})$ if its removal changes the optimal solution of $\mathrm{SP}(\mathcal{S})$.
\end{definition}

\begin{definition}[Helly's dimension]
Helly's dimension of the scenario problem $\mathrm{SP}(\mathcal{S})$ is the smallest integer $h$ that $h\geq_{\text{ess sup}\mathcal{S}\subseteq\Delta^N}|\mathcal{S(S)}|$ holds for any finite $N\geq1$, where $|\mathcal{S(S)}|$ is the number of support scenarios.
\end{definition}
The above definitions are used to quantify the number of samples that are needed to guarantee a risk level $\epsilon$ in CCO. Next, we make the following structural assumptions necessary to rule out the pathological cases of the optimization problem.

\begin{assumption}[Non-empty interior]
\label{assmpt:non-empty-interior}
For all sampled scenarios $\mathcal{S}$, the scenario problem $\mathrm{SP}(\mathcal{S})$ and robust scenario problems $\mathrm{RSP}(\mathcal{S})$ are feasible, and the feasibility region has a non-empty interior. Further, each such problem has a unique solution.
\end{assumption}
Assumption~\ref{assmpt:non-empty-interior} is applicable only for the cases when the problems considered are convex. Our assumptions are well-established in literature~\cite{campi2008exact} and ensure that the underlying optimization problems are well-behaved. Proposition~\ref{prop:violation-static} is a folklore result in scenario approach in a static environment with a convex feasible set and was first proven by~\cite{calafiore2006scenario}.

\begin{proposition}[Violation Probability~\cite{calafiore2006scenario}]
\label{prop:violation-static}
Under Assumption 1, the optimal solution to the scenario problem $\mathrm{SP}(\mathcal{S})$ satisfies:
\begin{equation*}
\label{firstN}
\mathbb{P}^N\left(\mathbb{V}\left(x_{\mathcal{N}}^*\right)>\epsilon\right) \leq \left(\begin{array}{c}
N \\
h
\end{array}\right) (1-\epsilon)^{N-h}\doteq\beta
\end{equation*}
where, the probability $\mathbb{P}^N$ is taken with respect to $\textit{N}^{th}$ scenario, and \textit{h} is the Helly's dimension of $\mathrm{SP}(\mathcal{S})$.
\end{proposition}

\begin{proposition}
The number of support scenarios for $\mathrm{SP}(\mathcal{S})$ is at most \textit{n} where \textit{n} is the dimension of $x$ and \textit{h} is Helly's dimension.
\end{proposition}

The rest of the paper will focus on robust scenario problems. The Helly's dimension for the \textit{robust scenario problem} $\text{RSP}(\mathcal{S})$ when the underlying constraint set is convex is upper bounded by $d$ (Theorem 3 in~\cite{calafiore2006scenario}). 
\section{Scenario Approach in Dynamic Environments}
\label{sec:nonstationary-scenario}
We consider the Robust Scenario Approach (\ref{RSP2}) wherein each scenario is sampled from a different and unknown measure. With this in mind, given a set of scenarios $\mathcal{S}:=\{\xi_1,\xi_2, \ldots, \xi_N\}, \ \xi_{i} \stackrel{iid}{\sim} \mathbb{P}_{i}$, we are interested in obtaining the sample complexity of ex-post violation probability under the measure of the (\textit{N}+1)-th scenario. To this end, we first extend Definition~\ref{def:violation-prob} to a dynamic environment when the scenarios are generated from time-varying distributions.
\begin{definition}[Violation in dynamic environments] 
The \textit{violation probability} of a candidate solution $x^*$ under uncertainty set with probability measure $\mathbb{P}_i$ is defined as the probability that $x^*$ is infeasible, i.e., $\mathbb{V}_i(x^*):=\mathbb{P}_i(x^*\notin \mathcal{X}_\xi)$
\end{definition}
In the rest of this paper, we focus on the relationship between the number of scenarios \textit {N} and the maximum value of \textit {violation probability} under a new scenario $\xi_{N+1} \sim \mathbb{P}_{N+1}$, i.e., $\mathbb{V}_{N+1} (x_{\mathcal{RS}}^*)\leq\epsilon$. where, $x_{\mathcal{RS}^\ast}$ is the solution to the robust scenario problem from scenarios $\xi_i \sim \mathbb{P}_i$. Our first step is to quantify the relationship between measures $\mathbb{P}_{i}, \mathbb{P}_{i+1}$. To this end, we will use the Wasserstein metric to couple the scenario-generating distributions.

\begin{definition} [Wasserstein  metric~\cite{villani2021topics}] 
The $1$-Wasserstein distance between two probability measures $\mathbb{P}$ and $\mathbb{P}^{\prime}$ defined on $\Xi$ is given by:
\begin{equation}
\begin{aligned}
d_{\mathrm{W}}\left(\mathbb{P}, \mathbb{P}^{\prime}\right)&:=\inf _{\pi \in \Pi\left(\mathbb{P}, \mathbb{P}^{\prime}\right)} \int_{\Xi \times \Xi}\left\|\xi-\xi^{\prime}\right\| \pi\left(d \xi, d \xi^{\prime}\right)\\
&=\inf _{\pi \in \Pi\left(\mathbb{P}, \mathbb{P}^{\prime}\right)} \underset{(\xi, \xi^{\prime}) \sim \pi}{\mathbb{E}}\left(\|x-y\|\right)
\end{aligned}
\end{equation}
where, $\|\cdot\|$ is a norm on $\Xi$, and $\Pi\left(\mathbb{P}, \mathbb{P}^{\prime}\right)$ denotes the set of all joint probability distributions of $\xi$ and $\xi^{\prime}$ with marginal distributions $\mathbb{P}$ and $\mathbb{P}^{\prime}$, respectively.
\end{definition}
In the rest of the paper, the $p$-Wasserstein distance denotes the Wasserstein distance under the norm $\|\cdot\|_p$. Based on the Wasserstein metric we propose a model for coupling different scenario-generating distributions. Our goal with such a model is to relate the scenario-generating probability measures across time.
\begin{definition}[Model A] 
We assume that the scenario-generating measures are constrained as follows:
\begin{equation}
\label{ModelA}
d_{\mathrm{W}}\left(\mathbb{P}_i, \mathbb{P}_j\right) \leq \rho
\end{equation}
for all pair of indices $(i,j), \ i\geq 1, j\geq1$. 
\end{definition}
A related but much more granular scenario-generating model is Model $B$. 
\begin{definition}[Model B] 
We assume
\begin{equation}
\label{ModelB}
d_{\mathrm{W}}\left(\mathbb{P}_i, \mathbb{P}_{j}\right) \leq \rho(i,j)
\end{equation}
for all indices $(i,j), \ i\geq 1, j \geq 1$, where $\rho: \mathbb{R}_+^2 \to \mathbb{R}_+$ is a known function satisfying $\rho(i,i)=0$.
\end{definition}
\textit{Model B} allows for a broader range of temporal shifts in the scenario-generating process, including gradual drifts over time, large but infrequent changes, or a combination thereof. All our results will be proven for scenario generation under Model $A$ or $B$. 

\subsection{Convex Scenario Problems}
We focus on the case when the underlying feasible set 
of the chance-constrained optimization is convex.
To this end, we make the following assumptions.  
\begin{assumption}[Convexity] 
The deterministic constraint $g(x) \leq 0$ is convex, and the random constraint $f(x,\xi)$ is convex in $x$ for every scenario of $\xi$. 
\end{assumption}

Our first result is a coupling-based inequality that relates to the probability that solutions remain feasible under consecutive scenarios.
\begin{lemma} 
\label{lem:coupling}
Let $x \in \chi$. For all $i \in [\textit{N}]$, for data generated under \textit{Model A}, we have:
\begin{equation}
\mathbb{P}_i(x_{\mathcal{RS}}^* \in \mathcal{X}_\xi)\leq\mathbb{P}_{N+1}(x_{\mathcal{RS}}^* \in \mathcal{X}_\xi)+\frac{\rho}{r_i}
\end{equation}
\end{lemma}
\begin{proof}
The proof of this lemma is based on~\cite{yan2022data}. Let $\left\{\nu_1, \nu_2, \cdots, \nu_N\right\}$ be a sequence of i.i.d. random variables such that: they are drawn independently from $\mathbb{P}_{N+1}(\cdot)$. and each pair of random variables $(\xi_i,\nu_i)$ are coupled to attain $1$-Wasserstein distance. Such a coupling exists by~\cite{wang2012coupling}. 
\begin{eqnarray}
\mathbb{P}_i(x_{\mathcal{RS}}^* \in \mathcal{X}_\xi)\stackrel{(a)}{=}&\mathbb{P}_i(x_{\mathcal{RS}}^* \in \mathcal{X}_\xi)\times(\,\Pr(\|\eta_i-\nu_i\| \leq r_i) \nonumber \\
&+\Pr(\|\eta_i-\nu_i\| > r_i))\, \nonumber \\
\leq&\mathbb{P}_i(x_{\mathcal{RS}}^* \in \mathcal{X}_\xi)\times \Pr(\|\eta_i-\nu_i\| \leq r_i) 
\nonumber \\
&+\Pr(\|\eta_i-\nu_i\| > r_i) \nonumber \\
\stackrel{(b)}{\leq}&\mathbb{P}_{N+1}(x_{\mathcal{RS}}^* \in \mathcal{X}_\xi)+\frac{\mathbb{E}(\|{\eta}_i-{\nu}_i\|)}{r_i} \nonumber \\
\leq&\mathbb{P}_{N+1}(x_{\mathcal{RS}}^* \in \mathcal{X}_\xi)+\frac{\rho}{r_i} \nonumber 
\end{eqnarray}
where, $(a)$ follows by definition of probability measure and $(b)$ follows since Markov's inequality. 
\end{proof}
When the scenario generating model is $B$, Corollary~\ref{cor:coupling} can be established using Lemma~\ref{lem:coupling}. 
\begin{corollary}
\label{cor:coupling}
Let $x \in \chi$. For all $i \in [\textit{N}]$, supposing data generated under \textit{Model B}, it holds that
\begin{equation}
\mathbb{P}_i(x_{\mathcal{RS}}^* \in \mathcal{X}_\xi)\leq\mathbb{P}_{N+1}(x_{\mathcal{RS}}^* \in \mathcal{X}_\xi)+\frac{\rho(i,N+1)}{r_i}
\end{equation}
\end{corollary}

Armed with Lemma~\ref{lem:coupling}, we are now ready to prove the violation risk when the scenarios are generated in a dynamic environment and the underlying problem is convex.

\begin{theorem} 
\label{thm:violation-convex}
Under Assumption 1 and 2, let $x_{\mathcal{RS}}^*$ be the optimal solution to the robust scenario approach $\mathrm{RSP}(\mathcal{S})$, supposing data generated under \textit{Model A}. Then:
\begin{equation}
\begin{aligned}
\underset{i=1,...,N}{\Pi}\mathbb{P}_i(\mathbb{V}_{N+1}(x_{\mathcal{RS}}^*)>\epsilon)<{N \choose d} \exp\Big(\frac{\rho(N+1-d)}{r_{\min}}-\epsilon\Big) \nonumber 
\end{aligned}
\end{equation}
\end{theorem}
\begin{proof} 
From Lemma~\ref{lem:coupling}, we have that:
\begin{equation}
\begin{aligned}
\mathbb{V}_i(x_{\mathcal{RS}}^*)&=1-\mathbb{P}_i(x_{\mathcal{RS}}^* \in \mathcal{X}_\xi)\\
&\geq 1-\mathbb{P}_{N+1}(x_{\mathcal{RS}}^* \in \mathcal{X}_\xi)-\frac{\rho}{r_i}\\
&>(\epsilon-\frac{\rho}{r_i})_+
\end{aligned}
\end{equation}
\indent Note that $\mathbb{V}_{N+1}(x_{\mathcal{RS}}^*)$ is a random variable since the solution $x_{\mathcal{RS}}^*$ of $\mathrm{RSP}(\mathcal{S})$ is, due to the fact that it depends on the random observation results $\eta_1, \eta_2, \cdots, \eta_N$. Thus, $\mathbb{V}_{N+1}(x_{\mathcal{RS}}^*)\leq\epsilon$ may hold for certain extractions $\eta_1, \eta_2, \cdots, \eta_N$, while $\mathbb{V}_{N+1}(x_{\mathcal{RS}}^*)>\epsilon$ may be true for others. The following proof aims to quantify the probability of "bad" extractions. Given \textit{N} scenarios $\eta_1, \eta_2, \cdots, \eta_N$ from observation, select a subset $I=\{i_1,...,i_d\}$ of $d$ indices from $\{1,...,N\}$ and let $x_{{\mathcal{RS}}_I}^*$ be the optimal solution of the program
\begin{equation}
\begin{aligned}
\mathrm{RSP}(\mathcal{S}): &\min _{x \in \chi} c^{\top} x\\
\text { s.t. } &x \in \cap_{j=1}^h \mathcal{X}_{\eta_{i_j}}^{r_{i_j}}
\end{aligned}
\end{equation}
Based on $x_{{\mathcal{RS}}_I}^*$ we next introduce a subset $\Delta_I^N$ of the set $\Delta^N$ defined as 
\begin{equation}
\begin{aligned}
\Delta_I^N\doteq \{(\eta_1,...,\eta_N):x_{{\mathcal{RS}}_I}^*=x_{\mathcal{RS}}^*\}
\end{aligned}
\end{equation}
\indent Let now $I$ range over the collection $\mathcal{I}$ of all possible choices of $d$ indexes from $\{1,...,N\}$, implying $\mathcal{I}$ contains:
$\left(\begin{array}{l}
N \\
d
\end{array}\right)$ sets, and we have
\begin{equation}
\Delta^N=\bigcup_{I \in \mathcal{I}} \Delta_I^N
\end{equation}
\indent Next, let $B\doteq\{(\eta_1,...,\eta_N):\mathbb{V}_{N+1}(x_{\mathcal{RS}}^*)>\epsilon\}$ and $B_I\doteq\{(\eta_1,...,\eta_N):\mathbb{V}_{N+1}(x_{{\mathcal{RS}}_I}^*)>\epsilon\}$. We have
\begin{equation}
\begin{aligned}
\label{B=}
B&=B\cap\Delta^N\\
&=B \cap\left(\cup_{I \in \mathcal{I}} \Delta_I^N\right)\\
&=\cup_{I \in \mathcal{I}} \left(B \cap \Delta_I^N\right)\\
&=\cup_{I \in \mathcal{I}} \left(B_I \cap \Delta_I^N\right)
\end{aligned}
\end{equation}
A bound for $\mathbb{P}_N\mathbb{P}_{N-1}...\mathbb{P}_2\mathbb{P}_1(B)$ is now obtained by bounding $\Pr(B_I \cap \Delta_I^N)$ and then summing over $I \in \mathcal{I}$. Fix any $I$, e.g., $I_0=\{1,...,d\}$. The set $B_{I_0}=B_{\{1,...,d\}}$ is in fact a cylinder with base in the Cartesian product of the first $d$ constraint domains, i.e. the condition $\mathbb{V}_{N+1}(x_{\mathcal{RS}}^*)>\epsilon$ only involves the first $d$ constraints. And constraints brought by $\eta_{d+1},...,\eta_N$ must be satisfied by $x_{{\mathcal{RS}}_I}^*$, otherwise, we would not have $x_{{\mathcal{RS}}_I}^*=x_{\mathcal{RS}}^*$. Thus, by the fact that the extractions are independent, we have:
\begin{equation}
\begin{aligned}
\scriptstyle {\underset{i\in\mathcal{I} \setminus I_0}{\Pi}\mathbb{P}_j\{(\eta_{h+1},...,\eta_N): (\overline\eta_1,...,\overline\eta_h,\eta_{h+1},...,\eta_N)\in B_{I_0} \cap \Delta_{I_0}^N\}}\\
< \underset{i\in\mathcal{I} \setminus I_0}{\Pi}(1-\mathbb{V}_i(x_{\mathcal{RS}}^*))
\end{aligned}
\end{equation}
where $(\overline\eta_1,...,\overline\eta_d) \in$ base of the cylinder. Integrating over the base of the cylinder $B_{I_0}$, we then obtain
\begin{equation}
\begin{aligned}
\underset{i\in\mathcal{I} \setminus I_0}{\Pi}\mathbb{P}_i\{B_{I_0} \cap \Delta_{I_0}^N\} <& \underset{i\in\mathcal{I} \setminus I_0}{\Pi}(1-\mathbb{V}_i(x_{\mathcal{RS}}^*))\\
&\times\underset{k \in I_0}{\Pi}\mathbb{P}_k(\text{base of } B_{I_0})\\
\leq& \underset{i\in\mathcal{I} \setminus I_0}{\Pi}(1-\mathbb{V}_i(x_{\mathcal{RS}}^*))\\
<&\underset{i\in\mathcal{I} \setminus I_0}{\Pi}(1-(\epsilon-\frac{\rho}{r_i})_+)
\end{aligned}
\end{equation}
From (\ref{B=}), we finally arrive to the desired bound for $\underset{j=1,...,N}{\Pi}\mathbb{P}_i(B)$ 
\begin{equation}
\begin{aligned}
\underset{i=1,...,N}{\Pi}\mathbb{P}_i(B)
&<\underset{I\in \mathcal{I}}{\Sigma}\{\underset{i\in\mathcal{I} \setminus I}{\Pi}(1-(\epsilon-\frac{\rho}{r_i})_+)\} \\
&< {N \choose d} \left(1-\epsilon +\frac{\rho}{r_{\min}}\right)^{N+1-d} \\
&\le {N \choose d} \exp\left(\Big(\frac{\rho(N+1-d)}{r_{\min}}-\epsilon\Big)\right)
\end{aligned}
\end{equation}
\end{proof}
Identical violation probability can be established for scenarios generated from Model $B$ (Corollary~\ref{cor:sample-convex-model}).
\begin{corollary} 
\label{cor:sample-convex-model}
Under Assumption 1 and 2, let $x_{\mathcal{RS}}^*$ be the optimal solution to the robust scenario approach $\mathrm{RSP}(\mathcal{N})$, supposing data generated under \textit{Model B}, it holds that:
\begin{equation*}
\begin{aligned}
\scriptstyle {\underset{i=1,...,N}{\Pi}\mathbb{P}_i(\mathbb{V}_{N+1}(x_{\mathcal{RS}}^*)>\epsilon)<\underset{I\in \mathcal{I}}{\Sigma}\{\underset{i\in\mathcal{I} \setminus I}{\Pi}(1-(\epsilon-\frac{\rho(i,N+1)}{r_i})_+)\}}
\end{aligned}
\end{equation*}
\end{corollary}

\subsection{Non-convex Scenario Problems}
In this section, we focus on the sample complexity of guaranteeing a risk level when the underlying scenario problems are non-convex. Establishing such guarantees is fundamentally different since the convexity of the feasible set ensures that the cardinality of support constraints can be determined apriori to sample scenarios. As indicated before, the cardinality of the \textit{robust scenario problem} when the underlying set is convex is at most $d$. However, this is not the case when the underlying constraint set is non-convex. Therefore, techniques that were used to establish Theorem~\ref{thm:violation-convex} are not applicable to non-convex problems. Approaches based on VC-dimension~\cite{vidyasagar1997theory} are very conservative and those using the convex-hull of the non-convex feasible region~\cite{vidyasagar2001randomized} are applicable to a very small class of problems. ~\cite{calafiore2012mixed,esfahani2014performance} consider the case of mixed-integer optimization and show that the number of support constraints can be determined for mixed-integer problems. In lieu of support constraints when the feasible set is convex, we define an invariant set which is the minimal set of support constraints determined \textit{a posteriori, after} solving the scenario problem~\eqref{eqn:robust-scenario-compact} based on the work of~\cite{campi2018general}. Our contribution over the work of~\cite{campi2018general}
lies in extending their proof when the scenarios are generated using time-varying distributions using coupling arguments outlined in Lemma~\ref{lem:coupling}. We use the cardinality of this invariant set as a measure of problem complexity for the \textit{robust scenario problem} when the constraint set is non-convex.

\begin{definition}[Invariant Constraints] 
\label{defn:support-constraints}
Given a set of scenarios $\mathcal{S} = \{\xi_{1}, \xi_{2}, \ldots, \xi_{N}\}$, the set of invariant constraints of cardinality $k$, $\mathcal{I} \subseteq \mathcal{S}$ such that $\text{RSP}(\mathcal{S}) = \text{RSP}(\mathcal{I})$.  
\end{definition}

The cardinality of the invariant set can be determined \textit{after} obtaining the set of scenarios $\mathcal{S}$ and our results will be based on the minimal size of such a set.

\begin{theorem}
\label{thm:non-convexity-complexity}
Let $\beta \in (0,1)$, and assume scenario generating Model A and Assumption~1 and $|\mathcal{I}|$ be the cardinality of set of invariant constraints. The following violation properties holds: \[
\prod_{i=1}^{N} \mathbb{P}_i(\mathbb{V}_{N+1}(\text{opt}_\mathbb{A}(\mathcal{S}))\leq\epsilon(|\mathcal{I}|))\geq1-\beta,
\] where the function $\epsilon(|\mathcal{I}|)$ is defined as
\begin{equation}
\begin{aligned}
\label{nonconvex}
&\epsilon(N)=1 \\
&\sum_{k=0}^{N-1} \{ \sum_{I_k\in \mathcal{I}}\{\prod_{i\in\mathcal{I} \setminus I_k}(1-(\epsilon(k)-\frac{\rho}{r_i})_{+})\}\}=\beta
\end{aligned}
\end{equation}
where $I_k$ means any subset in $\mathcal{I}$ with $k$ dimension.
\end{theorem}

\begin{proof}
Consider a family of functions $\mathcal{A}_m:\Delta^m \rightarrow \Theta$, indexed by the size $m$ of the sample. Similarly, the invariant constraints can be found as a function $\mathcal{B}_N:(\eta_1,...,\eta_N) \rightarrow \{i_1,...,i_k\},i_1<...<i_k$, where $(\eta_{i_1},...,\eta_{i_k})$ is an invariant set.\\
\indent Suppose $I_k$ be a selection of $k$ indexes ${i_1,...,i_k},i_1<...<i_k$ from ${1,...,N}$, and let $\theta_{I_k}=\mathcal{A}_k(\eta_{i_1},...,\eta_{i_k})$. Consider the subsets from a partition of $\Delta^N$ defined as follows:
\begin{equation}
\Delta^N_k=\{(\eta_1,...,\eta_N) \in \Delta^N :|\mathcal{B}_N(\eta_1,...,\eta_N)|=k\}
\end{equation}
We define the set $\Delta^N_{k,I_k}\subseteq\Delta^N_k$ by the following rule:$(\eta_1,...,\eta_N)\in \Delta^N_{k,I_k}$ if and only if $|\mathcal{B}_N(\eta_1,...,\eta_N)|=I_k$. Then we have
\begin{equation}
\Delta^N=\bigcup_{k=0}^{N}\bigcup_{I_k} \Delta^N_{k,I_k}
\end{equation}
Let $B=\{(\eta_1,...,\eta_N)\in\Delta^N:\mathbb{V}_{N+1}(\theta_N)>\epsilon(|\mathcal{I}|)\}$, and $B_{I_k}=\{(\eta_1,...,\eta_N)\in\Delta^N:\mathbb{V}_{N+1}(\theta_{I_k})>\epsilon(k)\}$. It holds that:
\begin{equation}
B=\Delta^N \cap B=\bigcup_{k=0}^{N-1}\bigcup_{I_k} \Delta^N_{k,I_k} \cap B_{I_k}
\end{equation}
\indent The above proof process is following \cite{campi2018general}, but the violation probability is different for different scenarios in a dynamic environment, i.e., $\mathbb{V}_{i}(\theta_{I_k})>\epsilon(k)-\frac{\rho}{r_i}$. Fixing any $I_k$, e.g.$I_k=\{1,...,k\}$, and supposing $(\bar{\eta}_1,...,\overline{\eta}_k)$ is the base of a cylinder.
\begin{equation}
\begin{aligned}
&\prod_{i=k+1}^N \mathbb{P}_i\left\{\left(\eta_{k+1}, \ldots, \eta_{N}\right):\left(\bar{\eta}_{1}, \ldots, \bar{\eta}_{k}, \eta_{k+1}, \ldots, \eta_{N}\right)\right. \\
&\left.\quad \in \Delta_{k,\{1, \ldots, k\}}^N \cap B_{\{1, \ldots, k\}}\right\} \\
&\leq \prod_{i=k+1}^N \mathbb{P}_i\left\{\bigcap_{i=k+1}^N\left\{\left(\eta_{k+1}, \ldots, \eta_{N}\right): \theta_{\{1, \ldots, k\}} \in \Theta_{\eta_{i}}\right\}\right\} \\
&=\prod_{i=k+1}^N \mathbb{P}_i\left\{\eta_{i}: \theta_{\{1, \ldots, k\}} \in \Theta_{\eta_{(i)}}\right\} \\
&\leq \prod_{i=k+1}^N(1-(\epsilon(k)-\frac{\rho}{r_i})_+)
\end{aligned}
\end{equation}
Integrating over the base of the cylinder $B_{\{1,..,k\}}$
\begin{equation}
\begin{aligned}
\prod_{i=1}^N \mathbb{P}_i &\{\Delta_{k,\{1, \ldots, k\}}^N \cap B_{\{1, \ldots, k\}}\}\\
&\leq \prod_{i=1}^k\{\text{base of} B_{\{1,..,k\}}\} \times \mathbb{P}_i\prod_{i=k+1}^N(1-(\epsilon(k)-\frac{\rho}{r_i})_+)\\
&\leq \prod_{i=k+1}^N(1-(\epsilon(k)-\frac{\rho}{r_i})_+)
\end{aligned}
\end{equation}
Considering that $I_k=\{1,...,k\}$ was made for only one fixed case, for any $I_k$, we obtain that
\begin{equation}
\begin{aligned}
&\prod_{i=1}^N \mathbb{P}_i\{\mathbb{V}_{N+1}(\theta_N)>\epsilon(|\mathcal{I}|)\}=\prod_{i=1}^N \mathbb{P}_i\{B\}\\
&\leq \sum_{k=0}^{N-1} \{ \sum_{I_k\in \mathcal{I}}\{\prod_{i\in\mathcal{I} \setminus I_k}(1-(\epsilon(|\mathcal{I}|)-\frac{\rho}{r_i})_+)\}\}=\beta
\end{aligned}
\end{equation}
\end{proof}

\begin{corollary}
Under the assumptions in Theorem~\ref{thm:non-convexity-complexity}, if the tolerable observation errors $r_i$ is set to a constant value $r_0$, it is easy to choose $\epsilon(k)$ by splitting $\beta$ evenly among the $N$ terms in the sum (\ref{nonconvex}) is 
\begin{equation}
\label{Choice}
\epsilon(k):= \begin{cases}1 & \text { if } k=N \\
1+\frac{\rho}{r_0}-\sqrt[N-k]{\frac{\beta}{N\left(\begin{array}{l}
{\scriptstyle N}\\
{\scriptstyle k}
\end{array}\right)}} & \text { otherwise }\end{cases}
\end{equation}    
\end{corollary} 

\indent \textit {\textbf{Remark 2} (Other Choices of $\epsilon(k)$)}: Notice that (\ref{Choice}) is only one simple choice of $\epsilon(k)$ under special case. Other choices than (\ref{Choice}) are possible and accessible by algorithms.

\section{Numerical Experiments}
\label{sec:case-study}
In this section, we provide simulation results for the sample complexity results established in previous sections for both convex and non-convex feasible regions. 
\subsection{Convex Feasible Region: The Probabilistic Point Covering Problem}
\begin{figure}[h] 
    \centering
    \includegraphics[width=0.4\textwidth]{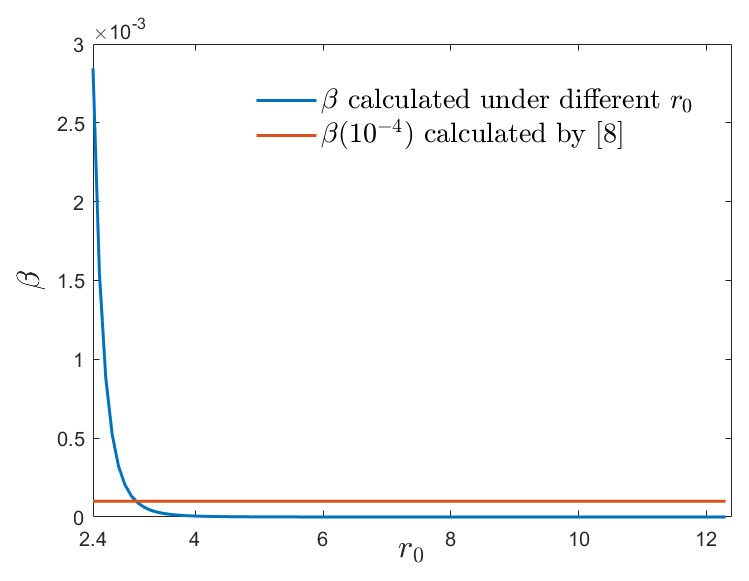}
    \caption{The relationship between the measurement error term and confidence parameter}
    \label{fig:beta}
\end{figure}

\indent Consider the one-dimensional probabilistic point covering problem, where the smallest interval is needed to cover the random variable with high probability. This problem can be written as a chance-constrained optimization form
\begin{equation}
\begin{array}{ll}
\min _{x, 
\gamma} & \gamma \\
\text { s.t. } & \mathbb{P}_{\xi}(x-\gamma \leq \xi \leq x+\gamma) \geq 1-\epsilon \\
& \gamma \geq 0
\end{array}
\end{equation}
where random samples of $\xi \in \mathbb{R}$ are supposed to be drawn from a fixed distribution, and the interval $[x-\gamma,x+\gamma]$ is what we seek for covering $1-\epsilon$ of the random variables $\xi$.
\indent The above problem can be solved by standard scenario approach (\ref{SP})~\cite{geng2019data}. In this section, we consider this point covering problems under time-varying distribution, and we need to guarantee our result satisfies the next generated ($N$+1)th data only from $N$ history scenarios, which can be formulated as 
\begin{equation}
\begin{array}{ll}
\label{PPC2}
\min _{x, \gamma} & \gamma \\
\text { s.t. } & \mathbb{P}_{N+1}(x-\gamma \leq \xi \leq x+\gamma) \geq 1-\epsilon \\
& \gamma \geq 0
\end{array}
\end{equation}
\indent We first solve the above problem by traditional scenario approach in~\cite{campi2008exact}. It is clear that the number of support scenarios is equal to the number of decision variables in (\ref{PPC2}), i.e. $h=2$. After setting risk level $\epsilon=0.1$, and confidence parameter $\beta=10^{-4}$, the number of needed scenarios $N=309$ is given by \textit{Theorem 1} in (\ref{firstN}). After setting the tolerable observation errors $r_i$, more robust results and smaller confidence parameters are generated by our new method to overcome the uncertainty in a dynamic environment. Table \ref{RR} compares the change in solution and the related confidence parameter $\beta$ calculated under different constant tolerable observation errors $r_0$ with the same tested scenarios.\\
\begin{table}[h] 
\caption{The optimal solution given by two methods}
\resizebox{\columnwidth}{!}{%
\begin{tabular}{ccccc}
\textbf{Measurement Error $(r)$} & $r_0=1.8$                                                   & $r_0=2$                                                      & $r_0=2.2$                                                  & $r_0=2.4$                                                   \\ \hline
\textbf{Scenario Approach in \cite{calafiore2006scenario}}           & \multicolumn{4}{c}{\begin{tabular}[c]{@{}c@{}}$\gamma^*=3.25$,\quad$\beta=10^{-4}$  \end{tabular}}                                                                                                                                                   \\ \hline
\textbf{Our Method}                  & \begin{tabular}[c]{@{}c@{}}$\gamma^*=5.05$\\ $\beta=0.17$\end{tabular} & \begin{tabular}[c]{@{}c@{}}$\gamma^*=5.25$\\ $\beta=0.042$\end{tabular} & \begin{tabular}[c]{@{}c@{}}$\gamma^*=5.45$\\ $\beta=0.011$\end{tabular} & \begin{tabular}[c]{@{}c@{}}$\gamma^*=5.65$\\ $\beta=0.0028$\end{tabular}
\end{tabular}%
}\label{RR}
\end{table}
\indent It can be seen from the above results that the larger measurement error term means the more conservative results with higher confidence parameter $\beta$. Fig.\ref{fig:beta} shows the relationship between the measurement error term and confidence parameter.
\subsection{Non-convex Feasible region: Control with Quantized Inputs}
\begin{figure}[h]
    \centering
    \includegraphics[width=0.4\textwidth]{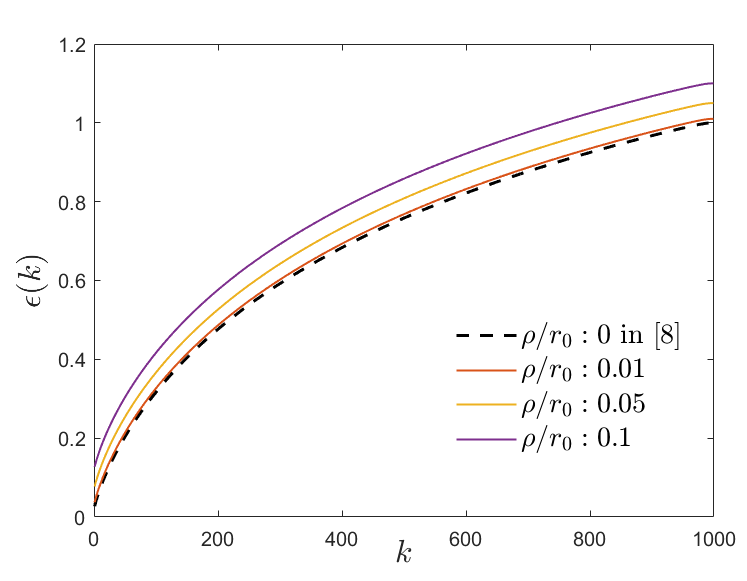}
    \caption{Comparing risk function between different $\rho$ and $r_0$ for mixed-integer optimal control}
    \label{fig:e_k}
\end{figure}

\indent We benchmark our results on a \textit{mixed-integer} optimal control problem introduced in \cite{campi2018general}, where the problem is the discrete-time control of an uncertain linear system.
\begin{equation}
x(t+1)=A x(t)+B u(t), \quad x(0)=x_0
\end{equation}
where the state variable is $x(t) \in \mathbb{R}^2$, with initial state $x_0=[1\quad1]^\top$, $B=[0\quad0.5]^\top$ is fixed, and uncertainty variable is $A \in \mathbb{R}^{2\times2}$, with independent Gaussian entries with standard deviation 0.02 each around the means
\begin{equation}
\bar{A}=\left[\begin{array}{cc}
0.8 & -1 \\
0 & -0.9
\end{array}\right]
\end{equation}
\indent Considering the actuation constraints, we suppose a finite input set: $u(t)\in \mathcal{U}:=\{-5,...,-1,0,1,...5\}$, which makes the control problem non-convex.
\indent Similar to \cite{campi2018general}, the control objective is to bring the state variable $x(t)$ to origin at time $T=8$ by choosing discrete $u(t)$ from $\mathcal{U}$. Since $x(T)= A^T x_0+\sum_{t=0}^{T-1} A^{T-1-t} B u(t)$, if we set 
\[
R=\left[\begin{array}{llll}
B & A B & \cdots & A^{T-1} B
\end{array}\right]
\]
and 
\[
\boldsymbol{u}=\left[\begin{array}{llll}
u(T-1) & u(T-2) & \cdots & u(0)
\end{array}\right]^{\top}
\]
The optimal control problem can then be reformulated as:
\begin{equation}
\begin{aligned}
\label{optimalcontrol}
\min _{h \in \mathbb{R}, \boldsymbol{u} \in \mathcal{U}} &h \\
\text { s.t. }&\left\|A^T x_0+R \boldsymbol{u}\right\|_{\infty} \leq h
\end{aligned}
\end{equation}
where $\| \cdot \|_\infty$ is the infinity norm. Because the matrices $A$ and $B$ both consist of uncertainty entries, the control design must incorporate some robustness. Under the scenario approach scheme, we treat these uncertainties by applying constraints to $N$ randomly extracted scenarios $A_{(i)}, i=1,...,N$, namely:
\begin{equation}
\begin{aligned}
\label{scenariocontrol}
\min _{h \in \mathbb{R}, \boldsymbol{u} \in \mathcal{U}} &h \\
\text { s.t. }&\left\|A_{(i)}^T x_0+R_{(i)} \boldsymbol{u}\right\|_{\infty} \leq h, \text{ for } i=1,...,N
\end{aligned}
\end{equation}
\indent We consider Model $A$ for generating the Gaussian entries in matrix $A$, where parameter $\rho$ is used to represent the shift in distribution.
\indent The scenario approach for the non-convex problem is a-posterior feasibility guarantees method, which means the risk level $\epsilon(|\mathcal{I}|)$ can only be given by Algorithm $\mathbb{B}$ after the solution under the given number of scenarios. If $N=1000$, and the tolerable observation error $r_i$ is set to a constant value $r_0$, then the relationships between $\epsilon(k)$ and $k$ under different parameter $\rho$ and $r_0$ are given in Fig.\ref{fig:e_k}.
\\
\indent Under the same invariant set cardinality $k$, it can be seen from Fig.\ref{fig:e_k} that the risk level $\epsilon(k)$ declines as the $\rho/r$ reduces, and when $\rho/r$ goes to 0, our results are identical with~\cite{campi2018general}. This numerically corroborates that our results extend the results of~\cite{campi2018general} to dynamic environments. 

\section{Conclusion}
In this paper, we propose novel sample complexity bounds on scenario generation in dynamic environments. To the best, this is the first work to provide bounds on the number of scenarios necessary to guarantee ex-post risk levels in chance-constrained optimization when the scenarios are generated from a dynamic environment. We provide guarantees for both cases when the underlying chance-constrained problem is convex or non-convex. We believe that this paper significantly contributes to antecedent literature on scenario generation and provides a holistic extension opening several avenues for future work. Of particular interest to us is the development of novel models for scenario-generating measures. Another possible avenue is the use of these guarantees in practical applications.

\bibliographystyle{IEEEtran}
\bibliography{ref}
\section*{Appendices}
\subsection{The Analytical Expression of 1-Wasserstein Distance Between Two Normal Distribution} \label{AppendicesA}
\begin{figure}[h]
    \centering
    \includegraphics[width=0.4\textwidth]{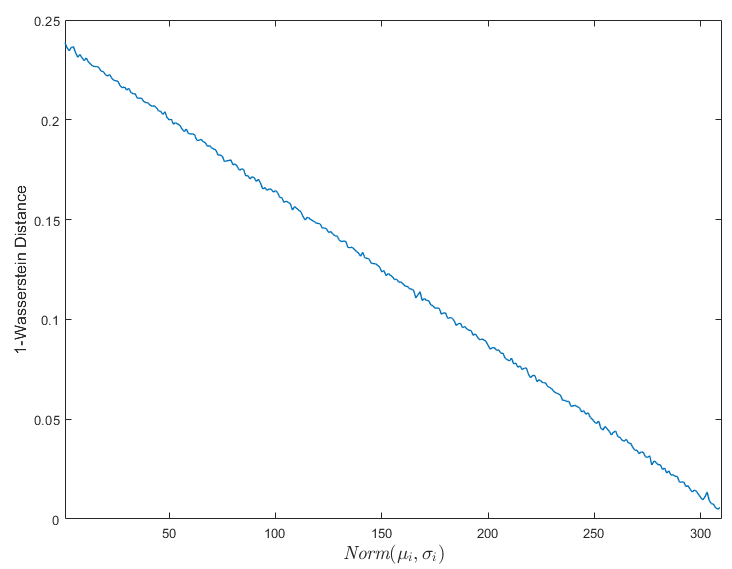}
    \caption{The approximation of 1-Wasserstein distance between $\textit{Norm}(\mu_i,\sigma_i)$ and $\textit{Norm}(\mu_{310},\sigma_{310})$}
    \label{fig:1-Wasserstein distance}
\end{figure}

We first introduce a lemma that characterizes the difference between two independent normal variables.
\begin{lemma}[\cite{weisstein2002normal}] Let $X$ and $Y$ are two independent normal variables with means and variances $(\mu_x,\sigma_x)$ and $(\mu_y,\sigma_y)$, respectively. Then $X-Y$ will follow a normal distribution with mean $\mu_x-\mu_y$ and variance $(\sigma_x^2,\sigma_y^2)$.
\end{lemma}
Supposing we want to calculate the 1-Wasserstein distance between one-dimensional normal distribution $\xi_i\sim\textit{Norm}(\mu_i,\sigma_i)$ and $\xi_{N+1}\sim\textit{Norm}(\mu_{N+1},\sigma_{N+1})$ under the shifting rules (\ref{rules}). Then we have
\begin{equation}
\begin{aligned}
d_{\mathrm{W}}\left(\mathbb{P}_i, \mathbb{P}_{N+1}\right)&=\inf _{\pi \in \Pi\left(\mathbb{P}_i, \mathbb{P}_{N+1}\right)} \underset{(\xi_i, \xi_{N+1}) \sim \pi}{\mathbb{E}}\left(|x-y|\right)\\
&=|\mu_i-\mu_{N+1}|\\
&=0.2\times(1-\frac{i-1}{N})
\end{aligned}
\end{equation}
Fixed the number of needed scenarios $N=309$, and based on the classic approximation algorithm presented in \cite{kantorovich2006translocation}, the approximation of 1-Wasserstein distance between $\textit{Norm}(\mu_i,\sigma_i)$ and $\textit{Norm}(\mu_{310},\sigma_{310})$ is shown in Fig.\ref{fig:1-Wasserstein distance}\\
\\
The introduced approximation method above is used to show the expansibility of our method in distribution-free circumstances.

\subsection{Simulation Details}
\begin{figure}[H]
    \centering
    \includegraphics[width=0.4\textwidth]{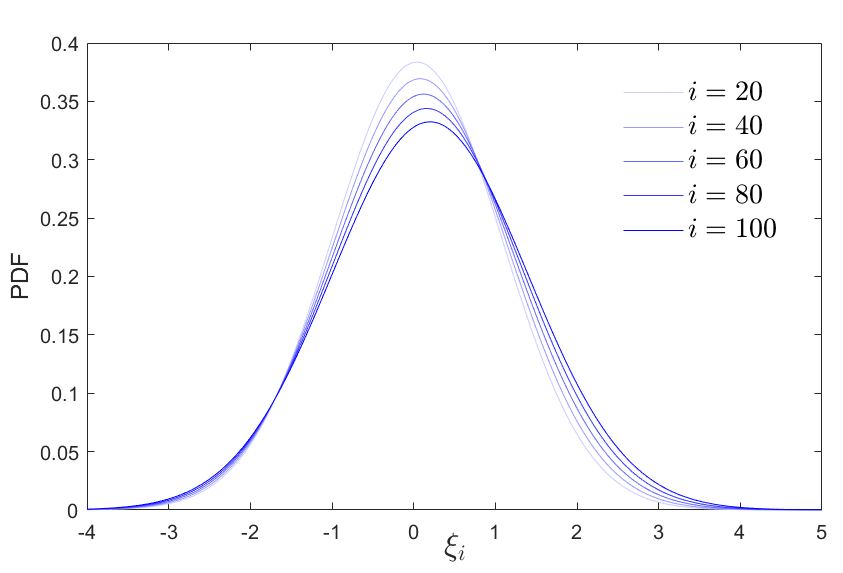}
    \caption{The changing distribution of $\xi_i$}
    \label{fig:datagene}
\end{figure}

\indent For better illustrating our ideas, we suppose the random variable $\xi_i\sim\textit{Norm}(\mu_i,\sigma_i)$ is distributed according to Normal distribution with both varying mean parameter and varying standard deviation following the rules below
\begin{equation}
\begin{aligned}
\label{rules}
&\mu_i=0.2 \times i/N \\
&\sigma_i=1+0.2 \times i/N 
\end{aligned}
\end{equation}
Supposing $N=100$, the probability density functions (PDF) of different $\xi_i$ are shown as Fig.\ref{fig:datagene}.\\
\indent Before applying the robust scenario approach scheme proposed in this paper, we first classify that the scenario generation process shown in Fig.\ref{fig:datagene} belongs to \textit{Model B} in (\ref{ModelB}), whose 1-Wasserstein distance between two one-dimensional distribution can be calculated by
\begin{equation}
\label{distance}
d_{\mathrm{W}}\left(\mathbb{P}_i, \mathbb{P}_{j}\right) = \int_{\mathbb{R}}\left|F_i(\xi)-F_{j}(\xi)\right| \mathrm{d}\xi
\end{equation}
where $F_i(\cdot)$ is the cumulative distribution functions (CDF) of Normal distribution $\textit{Norm}(\mu_i,\sigma_i)$.
\end{document}